\documentclass[preprints,article,accept,moreauthors]{my-mdpi} 

\firstpage{1} 
\pubvolume{7}
\issuenum{2}
\articlenumber{104}
\pubyear{2023}
\copyrightyear{2023}
\datereceived{30 November 2022} 
\daterevised{28 December 2022}
\dateaccepted{14 January 2023} 
\datepublished{18 January 2023} 
\hreflink{https://doi.org/\linebreak 10.3390/fractalfract7020104} 
\doinum{10.3390/fractalfract7020104}

\usepackage{amsmath,amssymb}
\usepackage{mathrsfs}

\allowdisplaybreaks 

\Title{Approximate Controllability of Delayed Fractional Stochastic Differential Systems with Mixed Noise and Impulsive Effects}

\TitleCitation{Approximate Controllability of Delayed Fractional Stochastic Differential Systems with Mixed Noise and Impulsive Effects}

\Author{Naima Hakkar $^{1}$\orcidA{}, 
Rajesh Dhayal $^{2}$\orcidB{}, 
Amar Debbouche $^{1,3,}$*\orcidC{} and 
Delfim F. M. Torres $^{3}$\orcidD{}}

\AuthorNames{Naima Hakkar, Rajesh Dhayal, Amar Debbouche and Delfim F. M. Torres}

\AuthorCitation{Hakkar, N.; Dhayal, R.; Debbouche, A.; Torres, D.F.M.}

\address{%
$^{1}$ \quad Department of Mathematics, Guelma University, Guelma 24000, Algeria; naima\_mt@yahoo.fr\\
$^{2}$ \quad School of Mathematics, Thapar Institute of Engineering and Technology, Patiala 147004, India; 
rajesh.dhayal@thapar.edu\\
$^{3}$ \quad CIDMA---Center for Research and Development in Mathematics and Applications, 
\mbox{Department of Mathematics}, University of Aveiro, 3810-193 Aveiro, Portugal; delfim@ua.pt}

\corres{Correspondence: amar\_debbouche@yahoo.fr}

\abstract{We herein report a new class of impulsive fractional stochastic differential 
systems driven by mixed fractional Brownian motions with infinite delay and Hurst 
parameter $\hat{\cal H} \in ( 1/2, 1)$. Using fixed point techniques, 
a $q$-resolvent family, and fractional calculus, we discuss the existence 
of a piecewise continuous mild solution for the proposed system. Moreover, 
under appropriate conditions, we investigate the approximate controllability 
of the considered system. Finally, the main results are demonstrated 
with an illustrative example.}



\begin{document}

\section{Introduction}

For a long time, the subject of fractional calculus and its 
applications has gained a lot of importance, mainly because fractional 
calculus has become a powerful tool with more accurate and successful results in
modeling several complex phenomena in numerous, seemingly diverse and widespread 
fields of science and engineering. It was found that various, especially 
interdisciplinary, applications can be elegantly modeled with the help of 
fractional derivatives \cite{F12,F11,F2,F1}. See also the recent 
works of \cite{F3,F4,F5,ref6}. 

Fractional Brownian motion (fBm for short) is a family of Gaussian random processes 
that are indexed by the Hurst parameter $\hat{\cal H}\in(0,1)$. It is a self-similar 
stochastic process with long-range dependence and stationary increment properties 
when $\hat{\cal H}>1/2$. For more recent works on fractional Brownian motion, 
see~\cite{S1,S2,S3,S4,S4.1,S5} and the references therein.

In order to describe various real-world problems in physical and engineering sciences 
subject to abrupt changes at certain instants during the evolution process, 
impulsive fractional differential equations have become important in recent years as
mathematical models of many phenomena in both physical and social sciences. Impulsive 
effects begin at any arbitrary fixed point and continue with a finite time interval, 
known as non-instantaneous impulses. For more details, we refer the reader to 
\cite{I1,I2,I3,I4,ref7,I5,I6,I7,I8}. 

\textls[25]{The concept of controllability plays a major role in finite dimensional control theory. 
However, its generalization to infinite dimensions is too strong and has limited applicability, 
while approximate controllability is a weaker concept completely adequate in applications \cite{Trig}.}

Recently, many authors have established approximate controllability results 
of (fractional) impulsive systems~\cite{A1,A2,A3,A5,A6,XDA,A8.1}. For example,
 Kumar et al.~\cite{A7} investigated the approximate controllability for impulsive 
semilinear control systems with delay; Anukiruthika et al.~\cite{A8} analyzed 
the approximate controllability of semilinear stochastic systems with impulses. 
Although several works exist in this area, the study of the approximate 
controllability of impulsive fractional stochastic differential systems 
driven by mixed noise with infinite delay and Hurst parameter $\hat{\cal H} \in ( 1/2, 1)$ 
is still an understudied topic in the literature. This fact provides 
the motivation of our current work.

We consider an impulsive fractional stochastic 
delay differential equation with mixed fractional Brownian motion
defined by
\begin{equation}
\begin{cases}
^cD_{t}^{q}z(t) = {\cal P}z(t)+{\cal F}(t,z_{t})
+{\cal G}(t,z_{t})\dfrac{d\hat{\cal W}(t)}{dt}
+\sigma(t)\dfrac{d{\cal B}^{\hat{\cal H}}(t)}{dt},
\quad t\in \cup_{i=0}^{m}(s_{i},t_{i+1}],& \text{}\\
z(t)  =  {\cal K}_{i}(t,z_{t}),
\quad t\in \cup_{i=1}^{m}(t_{i},s_{i}], 
&\text{}\\
z(t)= \phi(t),\;\;\phi(t)\in {\cal D}_{h}, 
& \text{}\label{eq1}
\end{cases}
\end{equation}
where ${\cal P}:{\cal D}({\cal P})\subset {\cal Z}\rightarrow {\cal Z}$ is the generator 
of an $q$-resolvent family $\{{\cal S}_{q}(t): t\geq 0\}$ on the separable Hilbert space 
${\cal Z}$, $^cD_{t}^{q}$ is the Caputo fractional derivative of order $1/2< q <1$, 
and state $z(\cdot)$ takes values in the space ${\cal Z}$, and 
$0= t_0 = s_0 < t_1 < s_1 < t_2<\cdots<t_m < s_m < t_{m+1} = T <\infty$. 
The functions ${\cal K}_{i}(t,z_{t})$ represent 
the non-instantaneous impulses during the intervals $(t_{i}, s_{i}],\,\,i = 1,2,\ldots,m$, 
$\hat{\cal W}=\{\hat{\cal W}(t):t\geq 0\}$ is a $Q$-Wiener process defined on a separable 
Hilbert space ${\cal Y}_{1}$, and ${\cal B}^{\hat{\cal H}}=\{{\cal B}^{\hat{\cal H}}(t):t\geq 0\}$ 
is a $Q$-fBm with the Hurst parameter $\hat{\cal H} \in (1/2, 1)$, defined on a separable 
Hilbert space ${\cal Y}_{2}$. The history-valued function $z_t:(-\infty,0] \rightarrow {\cal Z}$ 
is defined as $z_{t}(\theta)= z(t+\theta)$, $\forall \;\theta \leq 0$, and belongs to an abstract 
phase space ${\cal D}_{h}$. The initial data $\phi = \{\phi(t), t \in (\infty, 0]\}$ 
are ${\cal F}_{0}$-measurable, ${\cal D}_{h}$-valued random variable independent of 
$\hat{\cal W}$ and ${\cal B}^{\hat{\cal H}}$. The functions ${\cal F},\;{\cal G},$ $\sigma$, 
and ${\cal K}_{i}$ satisfy several suitable hypotheses, which will be specified later. 

The work is arranged as follows. In Section~\ref{sec:02}, relevant preliminaries 
are given that will be used later. In Section~\ref{sec:03}, we prove the existence 
of a piecewise continuous mild solution for the proposed system~(\ref{eq1}). 
Then, in Section~\ref{sec:04}, we study the approximate controllability 
for problem~(\ref{eq1}). In Section~\ref{sec:05}, 
an example is given to show the application of the obtained results. 
We end with Section~\ref{sec:06}, in which we present the conclusion of our results and also suggest
directions of possible future research. 

% ---------------------------------------------------------

\section{Preliminaries}
\label{sec:02}

Let $L({\cal Y}_{j},{\cal Z})$ denote the space of all linear and bounded operators 
from ${\cal Y}_{j}$ to ${\cal Z}$, $j=1,2$. The notation $\|\,\cdot\,\|$ represents 
the norms of ${\cal Z},$ ${\cal Y}_{j}$,  $L({\cal Y}_{j},{\cal Z})$. Let 
$(\Omega, {\cal F}, \{{\cal F}_{t}\}_{t \geq 0},P)$ be a filtered complete probability 
space, where ${\cal F}_{t}$ is the $\sigma$-algebra generated by 
$\{{\cal B}^{\hat{\cal H}}(e),\hat{\cal W}(e) : e \in [0,t]\}$ 
and $P$-null sets. Let ${\cal Q}_{j}\in L({\cal Y}_{j},{\cal Y}_{j})$ 
be the operators defined by ${\cal Q}_{j}e_{i}^{j} = \lambda_{i}^{j}e_{i}^{j}$ 
with finite trace $Tr(Q_{j}) = \sum_{i=1}^{\infty} \lambda_{i}^{j} < \infty$, 
where $\{\lambda_{i}^{j}\}_{i \geq 1}$ are non-negative real numbers 
and $\{e_{i}^{j}\}_{i\geq 1}$ is a complete orthonormal basis in ${\cal Y}_{j}$. 
Then, there exists  a real independent sequence $\mathscr{B}_{i}(t)$ 
of the standard Wiener process such that
$$
\hat{\cal W}(t)=\sum_{i=1}^{\infty}\sqrt{\lambda_{i}^{1}}\mathscr{B}_{i}(t)e_{i}^{1}.
$$
The infinite dimensional ${\cal Y}_{2}$-valued 
fBm ${\cal B}^{\hat{\cal H}}(t)$ is defined as 
$$
{\cal B}^{\hat{\cal H}}(t)=\sum_{i=1}^{\infty}
\sqrt{\lambda_{i}^{2}}\mathscr{B}^{\hat{\cal H}}_{i}(t)e_{i}^{2},
$$
where $\mathscr{B}^{\hat{\cal H}}_{i}(t)$ are real, independent fBms. 

Let $\mathscr{B}=\{\mathscr{B}(t), t \in {\cal J}\}$,
${\cal J}=[0,T]$ be a Wiener process and 
$\mathscr{B}^{\hat{\cal H}}=\{\mathscr{B}^{\hat{\cal H}}(t),
t \in {\cal J}\}$ be the one-dimensional fBm with Hurst index 
$\hat{\cal H}\in (1/2,1)$. The fBm $\mathscr{B}^{\hat{\cal H}}(t)$ has 
the following integral representation: 
$$
\mathscr{B}^{\hat{\cal H}}(t)
=\int_{0}^{t}\mathscr{K}_{\hat{\cal H}}(t,e)d\mathscr{B}(e),
$$
where the kernel $\mathscr{K}_{\hat{\cal H}}(t,e)$ is defined as
$$
\mathscr{K}_{\hat{\cal H}}(t,e)=\mathfrak{X}_{{\hat{\cal H}}}
e^{1/2-{\hat{\cal H}}}\int_{e}^{t}(\tau-e)^{{\hat{\cal H}}-3/2}
\tau^{{\hat{\cal H}}-1/2}d\tau\;\mbox{ for }\;t>e.
$$
We apply $\mathscr{K}_{\hat{\cal H}}(t,e) = 0$ if $t \leq e$. Note that
$\dfrac{\partial \mathscr{K}_{\hat{\cal H}}}{\partial t}(t,e)
=\mathfrak{X}_{{\hat{\cal H}}}\big(t/e\big)^{{\hat{\cal H}}-1/2}(t-e)^{{\hat{\cal H}}-3/2}$.
Here, $\mathfrak{X}_{{\hat{\cal H}}}=[{\hat{\cal H}}(2{\hat{\cal H}}-1)/\xi(2-2{\hat{\cal H}},
{\hat{\cal H}}-1/2)]^{1/2}$ and $\xi(\cdot,\cdot)$ is the Beta function.
For $\Lambda \in L^{2}([0,T])$, it follows from~\cite{bb} that the Wiener-type 
integral of the function $\Lambda$ w.r.t. fBm $\mathscr{B}^{\hat{\cal H}}$ is defined by
\begin{eqnarray*}
\int_{0}^{T} \Lambda(e)d\mathscr{B}^{\hat{\cal H}}(e)
=\int_{0}^{T} \mathscr{K}_{\hat{\cal H}}^{*}\Lambda(e)d\mathscr{B}(e),
\end{eqnarray*}
where $\mathscr{K}_{\hat{\cal H}}^{*}\Lambda(e)
=\int_{e}^{T}\Lambda(t)\dfrac{\partial \mathscr{K}_{\hat{\cal H}}}{\partial t}(t,e)dt$.

Let $\varphi_{j} \in L({\cal Y}_{j},{\cal Z})$ and define 
$$
\|\varphi_{j}\|_{{\cal L}^{j}_{2}}= \left[\sum_{i=1}^{\infty} 
\Big\|\sqrt{\lambda_{i}^{j}}\varphi_{j} e_{i}^{j}\|^{2}\right]^{1/2}.
$$
If $\|\varphi_{j}\|_{{\cal L}^{j}_{2}}< \infty$, then $\varphi_{j}$ are called 
$Q_{j}$-Hilbert--Schmidt operators, and the spaces ${\cal L}^{j}_{2}({\cal Y}_{j},{\cal Z})$ 
are real and separable Hilbert spaces with inner product 
$\langle \varphi^{1},\varphi^{2}\rangle_{{\cal L}^{j}_{2}}
= \sum_{i=1}^{\infty}\langle \varphi^{1} e_{i}^{j},\varphi^{2} e_{i}^{j}\rangle$. 
The stochastic integral of function $\Psi : {\cal J} 
\rightarrow {\cal L}^{2}_{2}({\cal Y}_{2},{\cal Z})$ 
w.r.t. fBm ${\cal B}^{\hat{\cal H}}$ is defined by
\begin{eqnarray}
\int_{0}^{t}\Psi(e)d{\cal B}^{\hat{\cal H}}(e)
=\sum_{i=1}^{\infty} \int_{0}^{t}\sqrt{\lambda_{i}^{2}}
\Psi(e)e_{i}^{2}d\mathscr{B}^{\hat{\cal H}}_{i}(e)
=\sum_{i=1}^{\infty} \int_{0}^{t}\sqrt{\lambda_{i}^{2}}
\mathscr{K}_{\hat{\cal H}}^{*}(\Psi e_{i}^{2})d\mathscr{B}_{i}(e).
\label{fbm}
\end{eqnarray}

\begin{Lemma}[See \cite{S1}] 
If $\Psi : {\cal J} \rightarrow {\cal L}^{2}_{2}({\cal Y}_{2},{\cal Z})$ 
satisfies $\int_{0}^{T}\|\Psi(e)\|_{{\cal L}^{2}_{2}}^{2}de< \infty$, 
then Equation~(\ref{fbm}) is a well-defined ${\cal Z}$-valued random variable 
such that
\begin{eqnarray*}
\mathbb{E}\bigg\|\int_{0}^{t}\Psi(e)d{\cal B}^{\hat{\cal H}}(e)\bigg\|^{2}
\leq 2{\hat{\cal H}}t^{2{\hat{\cal H}}-1}\int_{0}^{t}\|\Psi(e)\|_{{\cal L}^{2}_{2}}^{2}de.
\end{eqnarray*}
\end{Lemma}  

\begin{Lemma}[See \cite{S6.1}] 
For any $\alpha \geq 1$ and for an arbitrary ${\cal L}_{2}^{1}$-valued 
predictable process $\Upsilon(\cdot)$,  
\begin{eqnarray*}
\sup_{e\in[0,t]}\mathbb{E}\bigg\|\int_{0}^{e}
\Upsilon(\tau)d\hat{\cal W}(\tau)\Big\|^{2\alpha}
\leq (\alpha(2\alpha-1))^{\alpha}\bigg(\int_{0}^{t}\big(\mathbb{E}
\|\Upsilon(e)\|_{{\cal L}_{2}^{1}}^{2\alpha}\bigg)^{1/\alpha}de\bigg)^{\alpha},
\;t \in [0,T].
\end{eqnarray*} 
For $\alpha=1$, we obtain 
$$
\sup_{e\in[0,t]}\mathbb{E}\bigg\|\int_{0}^{e}\Upsilon(\tau)d\hat{\cal W}(\tau)
\bigg\|^{2}\leq \int_{0}^{t}\mathbb{E}\|\Upsilon(e)\|_{{\cal L}_{2}^{1}}^{2}de.
$$
\end{Lemma}  

Assume that $h : (-\infty, 0] \rightarrow (0, \infty)$ 
with $\varpi = \int_{-\infty}^{0} h(t)dt < \infty$ is a continuous 
function. We define ${\cal D}_{h}$ by
\vspace{-12pt} {}
\begin{adjustwidth}{-\extralength}{0cm}
\begin{eqnarray*}
{\cal D}_{h}&=&\Bigg\{\phi : (-\infty,0] \rightarrow {\cal Z},
\;\mbox{for\;\,any}\; a > 0,\;(\mathbb{E}|\phi(\theta)|^{2})^{1/2} \;
\mbox{{is a measurable and bounded function on}  
}\\
&&\;[-a, 0]\; \mbox{with}\; \phi(0) = 0,\;\mbox{and}\; 
\int_{-\infty}^{0}h(e)\sup_{e\leq\theta\leq 0}(\mathbb{E}|
\phi(\theta)|^{2})^{1/2}de < \infty \Bigg\}.
\end{eqnarray*}
\end{adjustwidth}
If ${\cal D}_{h}$ is endowed with the norm
\begin{eqnarray*}
\|\phi\|_{{\cal D}_{h}}= \int_{-\infty}^{0} h(e)
\sup_{e\leq\theta\leq 0}(\mathbb{E}\|\phi(\theta)\|^{2})^{1/2}de,
\;\phi\in {\cal D}_{h},
\end{eqnarray*}
then $({\cal D}_{h},\|\,.\,\|_{{\cal D}_{h}})$ is a Banach space~\cite{r7}.

Define the space ${\cal D}_{T}=\{z : (-\infty,T]\rightarrow {\cal Z}$, 
$z\vert_{{\cal J}_{i}}\in C({\cal J}_{i},{\cal Z})$, $i= 0,1,\ldots,m$, 
and there exist $z(t^{-}_i )$ and $z(t^{+}_i)$ with $z(t^{-}_i ) = z(t_{i}),$ 
and $z_{0}=\phi \in {\cal D}_{h}\},$ with the norm
\begin{eqnarray*}
\|z\|_{{\cal D}_{T}}= \|\phi\|_{{\cal D}_{h}} 
+ \sup_{t\in [0,T]}(\mathbb{E}\|z(t)\|^{2})^{1/2},
\end{eqnarray*}
where ${\cal J}_{i} = (t_{i}, t_{i+1}],\,i =0,1,\ldots,m$.

\begin{Lemma}[see \cite{A4}] 
\label{LM1}
If for all $t \in [0,T]$, $z_{t} \in {\cal D}_{h}$, 
$z_{0} \in {\cal D}_{h}$, then
\begin{eqnarray*}
\|z_{t}\|_{{\cal D}_{h}} \leq \varpi
\sup_{t\in [0,T]}(\mathbb{E}\|z(t)\|^{2})^{1/2}+\|z_{0}\|_{{\cal D}_{h}}.
\end{eqnarray*}
\end{Lemma}

\begin{Definition}[see \cite{TT}] 
Let ${\cal M} > 0$, $\theta \in [ \pi/2, \pi]$, 
and $\omega \in \mathbb{R}$. A closed and linear operator ${\cal P}$ 
is called a sectorial operator if 
\begin{enumerate}
\item $\rho({\cal P}) \subset \sum_{(\theta,\omega)} 
= \{\lambda \in \mathbb{C}: \lambda\neq\omega,|arg(\lambda-\omega)|<\theta\}$, 
\item $\|{\cal R}(\lambda,{\cal P})\| 
\leq {\cal M}/|\lambda-\omega|,\; \lambda \in \sum_{(\theta,\omega)}$.
\end{enumerate}
\end{Definition}

\begin{Lemma}[see \cite{DD1}] 
Let ${\cal P}$ be a sectorial operator.
Then, the unique solution of the linear fractional system
\begin{equation*}
\begin{cases}
^cD_{t}^{q}z(t)= {\cal P}z(t)+{\cal F}(t),
\quad t > t_{0} \geq 0,\quad 0<q<1, & \text{}\\
z(t)  = \phi(t),\quad t \leq t_{0}, & \text{}
\end{cases}
\end{equation*}
is given by
\begin{eqnarray*}
z(t)={\cal T}_{q}(t-t_{0})z(t_{0})
+ \int_{t_{0}}^{t} {\cal S}_{q}(t-e){\cal F}(e)de,
\end{eqnarray*}
where
\begin{eqnarray*}
{\cal T}_{q}(t) &=&\dfrac{1}{2\pi i}\int_{B_{r}}e^{\lambda t}
\dfrac{\lambda^{q-1}}{\lambda^{q}-{\cal P}}d\lambda,\\
{\cal S}_{q}(t)&=&\dfrac{1}{2\pi i}
\int_{B_{r}}\dfrac{e^{\lambda t}}{\lambda^{q}-{\cal P}}d\lambda.
\end{eqnarray*}
Here, $B_{r}$ denotes the Bromwich path.
\end{Lemma} 

% ---------------------------------------------------------

\section{Solvability Results} 
\label{sec:03}

We assume the following hypotheses.

\begin{Hypothesis}[\textbf{H1}] If $q \in (0,1)$ and ${\cal P} \in {\cal P}^{q}(\theta_{0}, \omega_{0})$, 
then, for any $z \in {\cal Z}$ and $t > 0$, we have $\|{\cal T}_{q}(t)\| 
\leq  C_{1}e^{\omega t}$ \mbox{and}  $\|{\cal S}_{q}(t)\| 
\leq C_{2}e^{\omega t}(1 + t^{q-1}),\;\omega > \omega_{0}$. Thus, we have
$$
\|{\cal T}_{q}(t)\| \leq  {\cal M}_{1}\;\;
\mbox{and}\;\; \|{\cal S}_{q}(t)\| \leq {\cal M}_{2}t^{q-1},
$$
where ${\cal M}_{1} = \sup_{0\leq t \leq T}C_{1}e^{\omega t}$ 
and ${\cal M}_{2} = \sup_{0\leq t \leq T} C_{2}e^{\omega t}(1 + t^{q-1})$. 
\end{Hypothesis} 

\begin{Hypothesis}[\textbf{H2}]There exists a constant $N_{{\cal F}}>0$ such that 
$$
\mathbb{E}\|{\cal F}(t,\psi_{1})-{\cal F}(t,\psi_{2})\|^{2} 
\leq N_{{\cal F}}\;\|\psi_{1}-\psi_{2}\|^{2}_{{\cal D}_{h}},
\quad \forall \;t \in {\cal J}, 
\quad \psi_{1},\psi_{2} \in {\cal D}_{h}.
$$
\end{Hypothesis} 

\begin{Hypothesis}[\textbf{H3}]Function $\sigma : {\cal J} 
\rightarrow {\cal L}^{2}_{2}({\cal Y}_{2},{\cal Z})$  
satisfies $\int_{0}^{t}\|\sigma(e)\|^{2}_{{\cal L}^{2}_{2}}de < \infty$, 
for every $t \in {\cal J}$, and there exists a constant 
$\Lambda_{\sigma}>0$ such that $\|\sigma(e)\|^{2}_{{\cal L}^{2}_{2}} 
\leq \Lambda_{\sigma}$, uniformly in ${\cal J}$. 
\end{Hypothesis}
\begin{Hypothesis}[\textbf{H4}]There exists a constant $N_{\cal G}>0$ such that
$$
\mathbb{E}\|{\cal G}(t,\psi_{1})-{\cal G}(t,\psi_{2})\|^{2}_{{\cal L}^{1}_{2}} 
\leq N_{\cal G}\;\|\psi_{1}-\psi_{2}\|^{2}_{{\cal D}_{h}},
\quad \forall \;t\in {\cal J}, 
\quad \psi_{1},\psi_{2} \in {\cal D}_{h}.
$$
\end{Hypothesis}
\begin{Hypothesis}[\textbf{H5}]There are constants $L_{{\cal K}_{i}}>0$, $i = 1,2,\ldots,m$, 
such that
\begin{eqnarray*}
\mathbb{E}\|{\cal K}_{i}(t,\psi_{1})-{\cal K}_{i}(t,\psi_{2}) \|^{2} 
\leq L_{{\cal K}_{i}}\;\|\psi_{1}-\psi_{2}\|^{2}_{{\cal D}_{h}},\,\,\,
\forall \,t \in {\cal J} ,\,\,\psi_{1},\psi_{2}  \in {\cal D}_{h}.
\end{eqnarray*}
\end{Hypothesis}

\begin{Definition} 
An ${\cal F}_{t}$-adapted random process $z:(-\infty,T] \rightarrow {\cal Z}$ 
is called the mild solution of~(\ref{eq1}) if, for every $t \in {\cal J}$, 
$z(t)$ satisfies $z_{0}=\phi \in {\cal D}_{h},$ $z(t)= {\cal K}_{i}(t,z_{t})$ 
for all $t \in (t_{i},s_{i}]$, $i=1,2,\ldots,m$, and 
\begin{eqnarray*}
z(t)&=&  \int_{0}^{t} {\cal S}_{q}(t-e){\cal F}(e,z_{e})de\\
&&+\int_{0}^{t}{\cal S}_{q}(t-e){\cal G}(e,z_{e})d\hat{\cal W}(e)
+ \int_{0}^{t}{\cal S}_{q}(t-e)\sigma(e)d{\cal B}^{\hat{\cal H}}(e),
\end{eqnarray*}
for all $t \in [0,t_{1}]$, and
\begin{eqnarray}
z(t)&=& {\cal T}_{q}(t-s_{i}){\cal K}_{i}(s_{i},z_{s_{i}}) 
+ \int_{s_{i}}^{t} {\cal S}_{q}(t-e){\cal F}(e,z_{e})de\nonumber\\
&&+\int_{s_{i}}^{t} {\cal S}_{q}(t-e){\cal G}(e,z_{e})d\hat{\cal W}(e)
+ \int_{s_{i}}^{t} {\cal S}_{q}(t-e)\sigma(e)d{\cal B}^{\hat{\cal H}}(e),
\label{mild1}
\end{eqnarray}
for all $t \in (s_{i},t_{i+1}]$, $i=1,2,\ldots,m$.
\end{Definition}

\begin{Theorem}
\label{TM11}
Assume that conditions {(H1)--(H5)} are satisfied. 
Then, problem~(\ref{eq1}) has a unique mild solution 
on $(-\infty ,T]$, provided that 
$$
L_{\cal R}= \max_{1\leq i \leq m}\left\{\eta_{0}, 
\varpi^{2}L_{{\cal K}_{i}}, \eta_{i}\right\}<1,
$$
where 
\begin{eqnarray*}
\eta_{0} & = &2{\cal M}_{2}^{2}\varpi^{2}\left(\dfrac{N_{\cal F}t_{1}^{2q}}{q^{2}}
+\dfrac{N_{\cal G}t_{1}^{2q-1}}{2q-1}\right),\\
\eta_{i} &=& \left(3{\cal M}_{1}^{2}L_{{\cal K}_{i}}
\varpi^{2}+3{\cal M}_{2}^{2}\varpi^{2}\left\{\dfrac{N_{\cal F}t_{i+1}^{2q}}{q^{2}}
+\dfrac{N_{\cal G}t_{i+1}^{2q-1}}{2q-1}\right\}\right).
\end{eqnarray*} 
\end{Theorem}

\begin{proof}
We define the operator $\Xi$ from 
${\cal D}_{T}$ to ${\cal D}_{T}$ as follows:
\begin{equation*}
(\Xi z)(t) =
\begin{cases}
\phi(t), & \text{} t \in (-\infty,0]\\
\int_{0}^{t} {\cal S}_{q}(t-e){\cal F}(e,z_{e})de\\
+\int_{0}^{t}{\cal S}_{q}(t-e){\cal G}(e,z_{e})d\hat{\cal W}(e)
+ \int_{0}^{t}{\cal S}_{q}(t-e)\sigma(e)d{\cal B}^{\hat{\cal H}}(e), 
& \text{} t \in [0,t_{1}]\\      
{\cal K}_{i}(t,z_{t}), & \text{} t \in (t_{i},s_{i}]\\    
{\cal T}_{q}(t-s_{i}){\cal K}_{i}(s_{i},z_{s_{i}}) 
+ \int_{s_{i}}^{t} {\cal S}_{q}(t-e){\cal F}(e,z_{e})de\\ 
+ \int_{s_{i}}^{t} {\cal S}_{q}(t-e){\cal G}(e,z_{e})d\hat {\cal W}(e)
+ \int_{s_{i}}^{t} {\cal S}_{q}(t-e)\sigma(e)d{\cal B}^{\hat{\cal H}}(e), 
& \text{} t \in (s_{i},t_{i+1}].
\end{cases}
\end{equation*}
For $\phi \in {\cal D}_{h}$, define
\begin{equation*}
g(t)=
\begin{cases}
\phi(t), & \text{}\ t \in (-\infty,0],\\
0, & \text{} \;\,t \in {\cal J}.
\end{cases}
\end{equation*}
Then, $g_{0}=\phi$. Next we define 
\begin{equation*}
\bar{y}(t)=
\begin{cases}
0, & \text{}\ t \in (-\infty,0]\\
y(t), & \text{} \;\,t \in {\cal J}
\end{cases}
\end{equation*}
for each $y \in C({\cal J},R)$ with $z(0) = 0$. 
If $z(\cdot)$ satisfies~(\ref{mild1}), then
$z(t) = g(t) + \bar{y}(t)$ for $t \in {\cal J}$, 
which implies that $z_{t} = g_{t} + \bar{y}_{t}$ 
for $t \in {\cal J}$, and the function $y(\cdot)$ 
satisfies
\begin{equation}
y(t)=
\begin{cases}
\int_{0}^{t} {\cal S}_{q}(t-e){\cal F}(e,g_{e} + \bar{y}_{e})de
+\int_{0}^{t}{\cal S}_{q}(t-e){\cal G}(e,g_{e} + \bar{y}_{e})
d\hat{\cal W}(e)\nonumber\\
+ \int_{0}^{t}{\cal S}_{q}(t-e)\sigma(e)d{\cal B}^{\hat{\cal H}}(e), 
& \text{} t \in [0,t_{1}]\\      
{\cal K}_{i}(t,g_{t} + \bar{y}_{t}), 
& \text{}t \in (t_{i},s_{i}]\\
{\cal T}_{q}(t-s_{i}){\cal K}_{i}(s_{i},g_{s_{i}}
+\bar{y}_{s_{i}}) + \int_{s_{i}}^{t} {\cal S}_{q}(t-e)
{\cal F}(e,g_{e} + \bar{y}_{e})de\nonumber\\
+\int_{s_{i}}^{t} {\cal S}_{q}(t-e){\cal G}(e,g_{e} 
+ \bar{y}_{e})d\hat{\cal W}(e)+ \int_{s_{i}}^{t} {\cal S}_{q}(t-e)
\sigma(e)d\hat{\cal B}^{\cal H}(e)
& \text{} t \in (s_{i},t_{i+1}].
\end{cases}
\end{equation}
Set ${\cal D}_{T}^{0}$ = $\big\{y \in {\cal D}_{T}$ 
such that $y_{0} = 0 \big\}$. For any $y \in {\cal D}_{T}^{0}$, 
we obtain
\begin{eqnarray*}
\|y\|_{{\cal D}_{T}^{0}}= \|y_{0}\|_{{\cal D}_{h}} 
+ \sup_{t\in {\cal J}}(\mathbb{E}\|y(t)\|^{2})^{1/2}
=\sup_{t\in {\cal J}}(\mathbb{E}\|y(t)\|^{2})^{1/2}.
\end{eqnarray*}
Thus, $({\cal D}_{T}^{0}, \|\,\cdot\,\|_{{\cal D}_{T}^{0}})$ is a Banach space.

Define the operator $\Psi$ from ${\cal D}_{T}^{0}$ 
to ${\cal D}_{T}^{0}$ as follows:
\begin{equation}
(\Psi y)(t)=
\begin{cases}
\int_{0}^{t} {\cal S}_{q}(t-e){\cal F}(e,g_{e} + \bar{y}_{e})de
+\int_{0}^{t}{\cal S}_{q}(t-e){\cal G}(e,g_{e} + \bar{y}_{e})
d\hat{\cal W}(e)\nonumber\\+ \int_{0}^{t}{\cal S}_{q}(t-e) 
\sigma(e)d{\cal B}^{\hat{\cal H}}(e), & \text{} t \in [0,t_{1}]\\        
{\cal K}_{i}(t,g_{t} + \bar{y}_{t}), & \text{} t \in (t_{i},s_{i}]\\
{\cal T}_{q}(t-s_{i}){\cal K}_{i}(s_{i},g_{s_{i}} + \bar{y}_{s_{i}}) 
+ \int_{s_{i}}^{t} {\cal S}_{q}(t-e){\cal F}(e,g_{e} + \bar{y}_{e})de\nonumber\\
+\int_{s_{i}}^{t} {\cal S}_{q}(t-e){\cal G}(e,g_{e} + \bar{y}_{e})d\hat{\cal W}(e)
+ \int_{s_{i}}^{t} {\cal S}_{q}(t-e)\sigma(e)d{\cal B}^{\hat{\cal H}}(e)
& \text{} t \in (s_{i},t_{i+1}].
\end{cases}
\end{equation}
In order to prove the existence result, we need to show that $\Psi$ 
has a unique fixed point. Let $y,y^{*} \in {\cal D}_{T}^{0}$. Then, 
for all $t \in [0, t_{1}]$, we have
\begin{eqnarray*}
\mathbb{E}\|({\Psi}y)(t)-(\Psi y^{*})(t)\|^{2} 
&\leq & 2\mathbb{E}\left\| \int_{0}^{t} {\cal S}_{q}(t-e)({\cal F}(e,g_{e} 
+ \bar{y}_{e})-{\cal F}(e,g_{e} + \bar{y}^{*}_{e}))de\right\|^{2}\qquad\qquad\quad\\
&&+ 2\mathbb{E}\left\|\int_{0}^{t} {\cal S}_{q}(t-e)({\cal G}(e,g_{e} 
+ \bar{y}_{e})-{\cal G}(e,g_{e} + \bar{y}^{*}_{e}))d\hat{\cal W}(e)\right\|^{2}\\
&\leq & \dfrac{2{\cal M}_{2}^{2}t_{1}^{q}}{q}\int_{0}^{t} (t-e)^{q-1}
N_{\cal F}\|\bar{y}_{e}-\bar{y}^{*}_{e}\|^{2}_{{\cal D}_{h}}de\\
&&+ 2{\cal M}_{2}^{2}\int_{0}^{t} (t-e)^{2q-2}N_{\cal G}\|\bar{y}_{e}
-\bar{y}^{*}_{e}\|^{2}_{{\cal D}_{h}}de\\
&\leq & \dfrac{2{\cal M}_{2}^{2}t_{1}^{q}}{q}\int_{0}^{t} (t-e)^{q-1}
N_{\cal F}\varpi^{2}\sup_{e \in {\cal J}}\mathbb{E}\|y(e)-y^{*}(e)\|^{2}de\\
&&+ 2{\cal M}_{2}^{2}\int_{0}^{t} (t-e)^{2q-2}N_{\cal G}\varpi^{2}
\sup_{e \in {\cal J}}\mathbb{E}\|y(e)-y^{*}(e)\|^{2}de\\
& \leq & 2{\cal M}_{2}^{2}\varpi^{2}\left(\dfrac{N_{\cal F}t_{1}^{2q}}{q^{2}}
+\dfrac{N_{\cal G}t_{1}^{2q-1}}{2q-1}\right)\|y-y^{*}\|^{2}_{{\cal D}_{T}^{0}}.
\end{eqnarray*}
Hence,
\begin{eqnarray}
\label{dd0}
\mathbb{E}\|({\Psi}y)(t)-(\Psi y^{*})(t)\|^{2} 
&\leq & 2{\cal M}_{2}^{2}\varpi^{2}\left(\dfrac{N_{\cal F}t_{1}^{2q}}{q^{2}}
+\dfrac{N_{\cal G}t_{1}^{2q-1}}{2q-1}\right)\|y-y^{*}\|^{2}_{{\cal D}_{T}^{0}}.
\end{eqnarray}

For $t \in (t_{i}, s_{i}]$, $i=1,2,\ldots,m,$ we have
\begin{eqnarray*}
\mathbb{E}\|({\Psi}y)(t)-(\Psi y^{*})(t)\|^{2} 
& \leq &\mathbb{E}\|{\cal K}_{i}(t,g_{t} + \bar{y}_{t})
-{\cal K}_{i}(t,g_{t} + \bar{y}^{*}_{t})\|^{2}\nonumber\\ 
&\leq & L_{{\cal K}_{i}} \|\bar{y}_{t}-\bar{y}^{*}_{t}\|^{2}_{{\cal D}_{h}}\nonumber\\ 
&\leq &  L_{{\cal K}_{i}}\varpi^{2}\sup_{t \in {\cal J}}\mathbb{E}\|y(t)-y^{*}(t)\|^{2}\nonumber\\
&\leq &  L_{{\cal K}_{i}}\varpi^{2}\|y-y^{*}\|^{2}_{{\cal D}_{T}^{0}}.\qquad\qquad\qquad\qquad
\end{eqnarray*}
Hence,
\begin{eqnarray}
\mathbb{E}\|({\Psi}y)(t)-(\Psi y^{*})(t)\|^{2} 
&\leq &   L_{{\cal K}_{i}}\varpi^{2}\|y-y^{*}\|^{2}_{{\cal D}_{T}^{0}}.
\end{eqnarray}
Similarly, for $t \in (s_{i},t_{i+1}]$, $i=1,2,\ldots,m$, we have
\begin{adjustwidth}{-\extralength}{0cm}
\begin{eqnarray*}
\mathbb{E}\|({\Psi}y)(t)-(\Psi y^{*})(t)\|^{2} 
&\leq & 3\mathbb{E}\|{\cal T}_{q}(t-s_{i})({\cal K}_{i}(s_{i},g_{s_{i}}
+\bar{y}_{s_{i}})-{\cal K}_{i}(s_{i},g_{s_{i}}+\bar{y}^{*}_{s_{i}}))\|^{2}
\qquad \qquad\qquad\\
&&+3\mathbb{E}\left\| \int_{s_{i}}^{t} {\cal S}_{q}(t-e)({\cal F}(e,g_{e} 
+ \bar{y}_{e})-{\cal F}(e,g_{e} + \bar{y}^{*}_{e}))de\right\|^{2}\\
&&+ 3\mathbb{E}\left\|\int_{s_{i}}^{t} {\cal S}_{q}(t-e)({\cal G}(e,g_{e} 
+ \bar{y}_{e})-{\cal G}(e,g_{e} + \bar{y}^{*}_{e}))d\hat{\cal W}(e)\right\|^{2}\\
&\leq & 3{\cal M}_{1}^{2}L_{{\cal K}_{i}}\varpi^{2}\|y-y^{*}\|^{2}_{{\cal D}_{T}^{0}}\\
&&+\dfrac{3{\cal M}_{2}^{2}t_{i+1}^{q}}{q}\int_{s_{i}}^{t} (t-e)^{q-1}
N_{\cal F}\|\bar{y}_{e}-\bar{y}^{*}_{e}\|^{2}_{{\cal D}_{h}}de\\
&&+ 3{\cal M}_{2}^{2}\int_{s_{i}}^{t} (t-e)^{2q-2}N_{\cal G}\|\bar{y}_{e}
-\bar{y}^{*}_{e}\|^{2}_{{\cal D}_{h}}de\\
&\leq & 3{\cal M}_{1}^{2}L_{{\cal K}_{i}}\varpi^{2}\|y-y^{*}\|^{2}_{{\cal D}_{T}^{0}}\\
&&+ \dfrac{3{\cal M}_{2}^{2}t_{i+1}^{q}}{q}\int_{s_{i}}^{t} (t-e)^{q-1}
N_{\cal F}\varpi^{2}\sup_{e \in {\cal J}}\mathbb{E}\|y(e)-y^{*}(e)\|^{2}de\\
&&+3{\cal M}_{2}^{2}\int_{s_{i}}^{t} (t-e)^{2(q-1)}N_{\cal G}\varpi^{2}
\sup_{e \in {\cal J}}\mathbb{E}\|y(e)-y^{*}(e)\|^{2}de\\
&\leq & \left(3{\cal M}_{1}^{2}L_{{\cal K}_{i}}\varpi^{2}
+3{\cal M}_{2}^{2}\varpi^{2}\left\{\dfrac{N_{\cal F}t_{i+1}^{2q}}{q^{2}}
+\dfrac{N_{\cal G}t_{i+1}^{2q-1}}{2q-1}\right\}\right)
\|y-y^{*}\|^{2}_{{\cal D}_{T}^{0}}.
\end{eqnarray*}
\end{adjustwidth}
Hence,

\begin{adjustwidth}{-\extralength}{0cm}
\begin{eqnarray}
\label{dd1}
\mathbb{E}\|({\Psi}y)(t)-(\Psi y^{*})(t)\|^{2} 
\leq \left(3{\cal M}_{1}^{2}L_{{\cal K}_{i}}\varpi^{2}
+3{\cal M}_{2}^{2}\varpi^{2}\left\{\dfrac{N_{\cal F}t_{i+1}^{2q}}{q^{2}}
+\dfrac{N_{\cal G}t_{i+1}^{2q-1}}{2q-1}\right\}\right)
\|y-y^{*}\|^{2}_{{\cal D}_{T}^{0}}.
\end{eqnarray}
\end{adjustwidth}
From Equations~(\ref{dd0})--(\ref{dd1}), we obtain that
\begin{eqnarray*}
\mathbb{E}\|{\Psi}y-\Psi y^{*}\|^{2}_{{\cal D}_{T}^{0}} 
\leq L_{\cal R}\|y-y^{*}\|^{2}_{{\cal D}_{T}^{0}},
\end{eqnarray*}
which implies that $\Psi$ is a contraction. Hence, $\Psi$ 
has a unique fixed point $y \in {\cal D}_{T}^{0}$, 
which is a mild solution of problem~(\ref{eq1}) 
on $(-\infty,T]$. 
\end{proof}

Next, using Krasnoselskii's fixed point theorem, we establish 
the second existence result. At this stage we make the following assumptions. 

\begin{Hypothesis}[\textbf{H6}]The map ${\cal F} : {\cal J} \times {\cal D}_{h} 
\rightarrow {\cal Z}$ is a continuous function, and there exists a 
continuous function $\xi_{1}:{\cal J} \rightarrow (0,\infty)$ such that
$$
\mathbb{E}\|{\cal F}(t,\psi)\|^{2} \leq \xi_{1}(t)\;\|\psi\|^{2}_{{\cal D}_{h}},
$$
for all \;$t\in {\cal J},$ and $\xi_{1}^{*}
= \sup_{t \in {\cal J}}\xi_{1}(t)$.
\end{Hypothesis}
\begin{Hypothesis}[\textbf{H7}]The map ${\cal G} : {\cal J} \times {\cal D}_{h} 
\rightarrow {\cal L}^{1}_{2}({\cal Y}_{1},{\cal Z})$ is a continuous 
function, and there exists a continuous function 
$\xi_{2}:{\cal J} \rightarrow (0,\infty)$ such that 
$$
\mathbb{E}\|{\cal G}(t,\psi)\|^{2}_{{\cal L}^{1}_{2}} 
\leq \xi_{2}(t)\;\|\psi\|^{2}_{{\cal D}_{h}},
$$
for all $t\in {\cal J}$ and 
$\xi_{2}^{*}= \sup_{t \in {\cal J}}\xi_{2}(t)$.
\end{Hypothesis}
\begin{Hypothesis}[\textbf{H8}]The inequality 
\begin{eqnarray*} 
L_{\cal HR}= 2{\cal M}_{2}^{2}\varpi^{2}\left(\dfrac{N_{\cal F}T^{2q}}{q^{2}}
+\dfrac{N_{\cal G}T^{2q-1}}{2q-1}\right) < 1
\end{eqnarray*}
holds and
$$
\max_{1\leq i \leq m}\left\{\kappa_{0}, \upsilon_{i}\lambda_{3}, \kappa_{i}\right\}<\pi,
$$
where 
\begin{eqnarray*}
\kappa_{0} & = & 3{\cal M}_{2}^{2}{t}_{1}^{2q}\left(\dfrac{\lambda_{1}}{q^{2}}
+\dfrac{\lambda_{2}}{t_{1}(2q-1)}+\dfrac{2\hat{\cal H}
\Lambda_{\sigma}t_{1}^{2\hat{\cal H}-2}}{2q-1}\right),\\
\kappa_{i} &=& 4{\cal M}_{1}^{2}\upsilon_{i}\lambda_{3}
+4{\cal M}_{2}^{2}{t}_{i+1}^{2q}\left(\dfrac{\lambda_{1}}{q^{2}}
+\dfrac{\lambda_{2}}{t_{i+1}(2q-1)}+\dfrac{2\hat{\cal H}\Lambda_{\sigma}
t_{i+1}^{2\hat{\cal H}-2}}{2q-1}\right).
\end{eqnarray*} 
\end{Hypothesis}
\begin{Hypothesis}[\textbf{H9}]The maps ${\cal K}_{i} : (t_{i},s_{i}] \times {\cal D}_{h} 
\rightarrow {\cal Z}$, $i=1,2,\ldots,m$, are continuous functions and 
\begin{enumerate}
\item[i.] there exist constants $\upsilon_{i}>0$, $i=1,2,\ldots,m$, such that 
$\mathbb{E}\|{\cal K}_{i}(t,\psi)\|^{2}
\leq\, \upsilon_{i}\,\|\psi\|^{2}_{{\cal D}_{h}}$
for all $t\in {\cal J}$;
\item[ii.] the set $\{b_{i}:b_{i}\in V(\pi,{\cal K}_{i})\}$ 
is an equicontinuous subset of $C((t_{i},s_{i}],{\cal Z})$,
$i=1,2,\ldots,m$, where $V(\pi,{\cal K}_{i})
=\{t\rightarrow {\cal K}_{i}(t,y_{t}): y\in  {\cal D}_{\pi}\}$.
\end{enumerate}

\end{Hypothesis}
The set ${\cal D}_{r}=\{y \in {\cal D}_{T}^{0}: \|y\|^{2}_{{\cal D}_{T}^{0}}
\leq r, r>0\}$ is clearly a  convex closed bounded set in ${\cal D}_{T}^{0}$ 
for each $y \in {\cal D}_{r}$. By Lemma~\ref{LM1}, we obtain
\begin{eqnarray*}
\|x_{t}+\bar{y}_{t}\|_{{\cal D}_{h}}^{2} 
&\leq & 2(\|x_{t}\|_{{\cal D}_{h}}^{2}+\|\bar{y}_{t}\|_{{\cal D}_{h}}^{2})\\
&\leq & 4\Big(\varpi^{2}\sup_{\nu\in [0,t]}\mathbb{E}\|x(\nu)\|^{2}
+\|x_{0}\|_{{\cal D}_{h}}^{2}\Big)+4\Big(\varpi^{2}\sup_{\nu\in [0,t]}
\mathbb{E}\|\bar{y}(\nu)\|^{2}+\|\bar{y}_{0}\|_{{\cal D}_{h}}^{2}\Big)\\
&\leq & 8(\|\phi\|_{{\cal D}_{h}}^{2}+\varpi^{2}r).
\end{eqnarray*}
Let
$$
\lambda_{1} = 8\xi_{1}^{*}(\|\phi\|_{{\cal D}_{h}}^{2}+\varpi^{2}r),\;\; 
\lambda_{2} = 8\xi_{2}^{*}(\|\phi\|_{{\cal D}_{h}}^{2}+\varpi^{2}r),\;\;
\lambda_{3} = 8(\|\phi\|_{{\cal D}_{h}}^{2}+\varpi^{2}r).
$$

\begin{Theorem}
\label{xxx}
Assume conditions {(H1)--(H9)} are satisfied. 
Then, problem~(\ref{eq1}) has at least one mild 
solution on $(-\infty, T]$.
\end{Theorem}

\begin{proof}
Let ${\cal E}_{1}$ : ${\cal D}_{r}$ $\rightarrow$ ${\cal D}_{r}$ 
and ${\cal E}_{2}$ : ${\cal D}_{r}$ $\rightarrow$ ${\cal D}_{r}$ 
be defined as
\begin{equation*}
{\cal E}_{1}(y)(t) =
\begin{cases}
0  & \text{} t \in [0,t_{1}]\\        
{\cal K}_{i}(t,g_{t} + \bar{y}_{t}), & \text{} t \in (t_{i},s_{i}]\\
{\cal T}_{q}(t-s_{i}){\cal K}_{i}(s_{i},g_{s_{i}} + \bar{y}_{s_{i}})
\qquad\qquad\qquad\qquad\qquad\qquad\qquad\qquad 
& \text{} t \in (s_{i},t_{i+1}]
\end{cases}
\end{equation*}
and
\begin{equation*}
{\cal E}_{2}(y)(t)=
\begin{cases}
\int_{0}^{t} {\cal S}_{q}(t-e){\cal F}(e,g_{e} + \bar{y}_{e})de\\
+\int_{0}^{t}{\cal S}_{q}(t-e){\cal G}(e,g_{e} + \bar{y}_{e})d\hat{\cal W}(e)
+ \int_{0}^{t}{\cal S}_{q}(t-e)\sigma(e)d{\cal B}^{\hat{\cal H}}(e), 
& \text{} \;\,t \in [0,t_{1}]\\    
0, & \text{}\ t \in (t_{i},s_{i}]\\                
\int_{s_{i}}^{t} {\cal S}_{q}(t-e){\cal F}(e,g_{e} + \bar{y}_{e})de\\
+\int_{s_{i}}^{t} {\cal S}_{q}(t-e){\cal G}(e,g_{e} + \bar{y}_{e})d\hat{\cal W}(e)
+\int_{s_{i}}^{t} {\cal S}_{q}(t-e)\sigma(e)d{\cal B}^{\hat{\cal H}}(e)
& \text{} \;\,t \in (s_{i},t_{i+1}].
\end{cases}
\end{equation*}
For convenience, we divide the proof into various steps.

Step 1. We show that ${\cal E}_{1}y + {\cal E}_{2}y^{*} \in {\cal D}_{r}$. 
For $y,y^{*} \in {\cal D}_{r}$ and for $t \in [0,t_{1}]$, we obtain
\vspace{-12pt}
\begin{adjustwidth}{-\extralength}{0cm}
\begin{eqnarray*}
\mathbb{E}\|({\cal E}_{1}y)(t) + ({\cal E}_{2}y^{*})(t)\|^{2} 
&\leq & 3\mathbb{E}\left\|\int_{0}^{t} {\cal S}_{q}(t-e){\cal F}(e,g_{e} 
+ \bar{y}^{*}_{e})de\right\|^{2}\\
&&+3\mathbb{E}\left\|\int_{0}^{t}{\cal S}_{q}(t-e){\cal G}(e,g_{e} 
+ \bar{y}^{*}_{e})d\hat{\cal W}(e)\right\|^{2}\\
&&+3\mathbb{E}\left\|\int_{0}^{t}{\cal S}_{q}(t-e)
\sigma(e)d{\cal B}^{\hat{\cal H}}(e)\right\|^{2}\\
&\leq & 3{\cal M}_{2}^{2}\left(\int_{0}^{t} (t-e)^{q-1}de\right)
\left(\int_{0}^{t} (t-e)^{q-1}\xi_{1}(e)\|g_{e} 
+ \bar{y}^{*}_{e}\|_{{\cal D}_{h}}^{2}de\right)\\
&& +3{\cal M}_{2}^{2}\int_{0}^{t} (t-e)^{2q-2}\xi_{2}(e)\|g_{e} 
+ \bar{y}^{*}_{e}\|_{{{\cal D}_{h}}}^{2}de\\
&& +6\hat{\cal H}\Lambda_{\sigma}{\cal M}_{2}^{2}t_{1}^{2
\hat{\cal H}-1}\int_{0}^{t} (t-e)^{2q-2}de\\
&\leq & 3{\cal M}_{2}^{2}{t}_{1}^{2q}\left(\dfrac{\lambda_{1}}{q^{2}}
+\dfrac{\lambda_{2}}{t_{1}(2q-1)}+\dfrac{2\hat{\cal H}\Lambda_{\sigma}
t_{1}^{2\hat{\cal H}-2}}{2q-1}\right).
\end{eqnarray*}
\end{adjustwidth}
Hence,
\begin{eqnarray}
\mathbb{E}\|({\cal E}_{1}y)(t) + ({\cal E}_{2}y^{*})(t)\|^{2} 
&\leq & 3{\cal M}_{2}^{2}{t}_{1}^{2q}\left(\dfrac{\lambda_{1}}{q^{2}}
+\dfrac{\lambda_{2}}{t_{1}(2q-1)}+\dfrac{2\hat{\cal H}
\Lambda_{\sigma}t_{1}^{2\hat{\cal H}-2}}{2q-1}\right).
\qquad\label{xx}
\end{eqnarray}
For $t \in (t_{i}, s_{i}]$, $i=1,2,\ldots,m$, we have
\begin{eqnarray*}
\mathbb{E}\|({\cal E}_{1}y)(t) + ({\cal E}_{2}y^{*})(t)\|^{2} 
& \leq & \mathbb{E}\|{\cal K}_{i}(t,g_{t} + \bar{y}_{t})\|^{2}\\
&\leq &  \upsilon_{i}\,\|g_{t} + \bar{y}_{t}\|^{2}_{{\cal D}_{h}}\\
&\leq & \upsilon_{i}\lambda_{3}.
\end{eqnarray*}
Hence,
\begin{eqnarray}
\mathbb{E}\|({\cal E}_{1}y)(t) + ({\cal E}_{2}y^{*})(t)\|^{2} 
& \leq &  \upsilon_{i}\lambda_{3}.
\end{eqnarray}
Similarly, for $t \in (s_{i}, t_{i+1}]$, $i=1,2,\ldots,m$, we have\vspace{-12pt}
\begin{adjustwidth}{-\extralength}{0cm}
\begin{eqnarray*}
\mathbb{E}\|({\cal E}_{1}y)(t) + ({\cal E}_{2}y^{*})(t)\|^{2} 
& \leq & 4\mathbb{E}\|{\cal T}_{q}(t-s_{i}){\cal K}_{i}(s_{i},g_{s_{i}} 
+ \bar{y}_{s_{i}})\|^{2}\\
&&+4\mathbb{E}\left\|\int_{s_{i}}^{t} {\cal S}_{q}(t-e){\cal F}(e,g_{e} 
+ \bar{y}^{*}_{e})de\right\|^{2}\\
&&+4\mathbb{E}\left\|\int_{s_{i}}^{t}{\cal S}_{q}(t-e){\cal G}(e,g_{e} 
+ \bar{y}^{*}_{e})d\hat{\cal W}(e)\right\|^{2}\\
&&+4\mathbb{E}\left\|\int_{s_{i}}^{t}{\cal S}_{q}(t-e)
\sigma(e)d{\cal B}^{\hat{\cal H}}(e)\right\|^{2}\\
&\leq & 4{\cal M}_{1}^{2}\upsilon_{i}\lambda_{3}\\
&&+4{\cal M}_{2}^{2}\left(\int_{s_{i}}^{t} (t-e)^{q-1}de\right)
\left(\int_{s_{i}}^{t} (t-e)^{q-1}\xi_{1}(e)\|g_{e} 
+ \bar{y}^{*}_{e}\|_{{\cal D}_{h}}^{2}de\right)\\
&& +4{\cal M}_{2}^{2}\int_{s_{i}}^{t} (t-e)^{2q-2}\xi_{2}(e)\|g_{e} 
+ \bar{y}^{*}_{e}\|_{{{\cal D}_{h}}}^{2}de\\
&&+8\hat{\cal H}\Lambda_{\sigma}{\cal M}_{2}^{2}
t_{i+1}^{2\hat{\cal H}-1}\int_{s_{i}}^{t} (t-e)^{2q-2}de\\
&\leq & 4{\cal M}_{1}^{2}\upsilon_{i}\lambda_{3}
+4{\cal M}_{2}^{2}{t}_{i+1}^{2q}\left(\dfrac{\lambda_{1}}{q^{2}}
+\dfrac{\lambda_{2}}{t_{i+1}(2q-1)}+\dfrac{2\hat{\cal H}
\Lambda_{\sigma}t_{i+1}^{2\hat{\cal H}-2}}{2q-1}\right).
\end{eqnarray*}
\end{adjustwidth}
Therefore,\vspace{-12pt}
\begin{adjustwidth}{-\extralength}{0cm}
\begin{eqnarray}
\mathbb{E}\|({\cal E}_{1}y)(t) + ({\cal E}_{2}y^{*})(t)\|^{2} 
&\leq & 4{\cal M}_{1}^{2}\upsilon_{i}\lambda_{3}+4{\cal M}_{2}^{2}{t}_{i+1}^{2q}
\left(\dfrac{\lambda_{1}}{q^{2}}+\dfrac{\lambda_{2}}{t_{i+1}(2q-1)}
+\dfrac{2\hat{\cal H}\Lambda_{\sigma}t_{i+1}^{2\hat{\cal H}-2}}{2q-1}\right).
\label{xx1}
\end{eqnarray}
\end{adjustwidth}
Equations~(\ref{xx})--(\ref{xx1}) imply that
\begin{eqnarray*}
\|{\cal E}_{1}y + {\cal E}_{2}y^{*}\|^{2}_{{\cal D}_{T}^{0}} 
&\leq & r.
\end{eqnarray*}
Thus, ${\cal E}_{1}y + {\cal E}_{2}y^{*} \in {\cal D}_{r}$.

Step 2. We show that the operator ${\cal E}_{1}$ 
is continuous on ${\cal D}_{r}$.
Let $\{y^{n}\}_{n=1}^{\infty}$ be a sequence such that $y^{n} \rightarrow y$ 
in ${\cal D}_{r}$. For all $t \in (t_{i}, s_{i}]$, $i=1,2,\ldots,m$, we have
\begin{eqnarray*}
\mathbb{E}\|({\cal E}_{1}y^{n})(t) - ({\cal E}_{1}y)(t)\|^{2}  
\leq  \mathbb{E}\| {\cal K}_{i}(t,g_{t} + \bar{y^{n}}_{t})
- {\cal K}_{i}(t,g_{t} + \bar{y}_{t})\|^{2}.
\end{eqnarray*}
Since the maps ${\cal K}_{i}$, $i=1,2,\ldots,m$, 
are continuous functions, one has 
\begin{eqnarray}
\lim_{n \rightarrow \infty} \|{\cal E}_{1}y^{n} 
- {\cal E}_{1}y\|^{2}_{{\cal D}_{T}^{0}}=0. 
\label{yy}
\end{eqnarray}
For all $t \in (s_{i}, t_{i+1}]$, $i=1,2,\ldots,m$, we have
\begin{eqnarray*}
\mathbb{E}\|({\cal E}_{1}y^{n})(t) - ({\cal E}_{1}y)(t)\|^{2} 
& \leq & \mathbb{E}\|{\cal T}_{q}(t-s_{i})({\cal K}_{i}(s_{i},
g_{s_{i}}+\bar{y^{n}}_{s_{i}})-{\cal K}_{i}(s_{i},g_{s_{i}}
+\bar{y}_{s_{i}})\|^{2}.
\end{eqnarray*}
Therefore, 
\begin{eqnarray}
\lim_{n \rightarrow \infty}\|{\cal E}_{1}y^{n} 
- {\cal E}_{1}y\|^{2}_{{\cal D}_{T}^{0}}=0. 
\label{yy1}
\end{eqnarray}
Equations~(\ref{yy}) and (\ref{yy1}) imply that the operator ${\cal E}_{1}$ 
is continuous on ${\cal D}_{r}$.

Step 3. The operator ${\cal E}_{1}$ maps bounded sets 
into bounded sets in ${\cal D}_{r}$. Let us show that 
for $r > 0$ there exists a $r > 0$ such that, for each 
$y \in {\cal D}_{r}$,  we obtain $\mathbb{E}\|{\cal E}_{1}(y)(t)\|^{2}$ 
$\leq r,$  for all $t \in (s_{i}, t_{i+1}],\;i=1,2,\ldots,m$.
For all $t \in (s_{i}, t_{i+1}]$, $i = 1,2,\ldots,m$, we have 
\begin{eqnarray*}
\mathbb{E}\|({\cal E}_{1}y)(t)\|^{2} 
& \leq & \mathbb{E}\|{\cal T}_{q}(t-s_{i}){\cal K}_{i}(s_{i},g_{s_{i}}
+\bar{y}_{s_{i}})\|^{2}\leq {\cal M}_{1}^{2}\upsilon_{i}\lambda_{3}.
\label{pp}
\end{eqnarray*}
For all $t \in (t_{i}, s_{i}]$, $i = 1,2,\ldots,m$, we have
\begin{eqnarray*}
\mathbb{E}\|({\cal E}_{1}y)(t)\|^{2} 
& \leq & \mathbb{E}\|{\cal K}_{i}(t,g_{t}+\bar{y}_{t})\|^{2}
\leq \nu_{i}\lambda_{3}.
\label{pp1}
\end{eqnarray*}
From the above equations, we obtain
\begin{eqnarray*}
\|{\cal E}_{1}y\|^{2}_{{\cal D}_{T}^{0}} & \leq & r,
\end{eqnarray*}
where $r=\max\{{\cal M}_{1}^{2}\upsilon_{i}\lambda_{3},\upsilon_{i}\lambda_{3}\}$. 
Hence, the operator ${\cal E}_{1}$ maps bounded sets into bounded sets 
in ${\cal D}_{r}.$

Step 4. The operator ${\cal E}_{1}$ is equicontinuous.
For all $\Delta_{1},\Delta_{2} \in (t_{i}, s_{i}]$, $\Delta_{1} < \Delta_{2}$, 
and $y \in {\cal D}_{r}$, we obtain
\begin{eqnarray}
\mathbb{E}\|({\cal E}_{1}y)(\Delta_{2}) - ({\cal E}_{1}y)(\Delta_{1})\|^{2}  
\leq \mathbb{E}\|{\cal K}_{i}(\Delta_{2},g_{\Delta_{2}} 
+ \bar{y}_{\Delta_{2}})-{\cal K}_{i}(\Delta_{1},
g_{\Delta_{1}} + \bar{y}_{\Delta_{1}})\|^{2}.
\label{zz}
\end{eqnarray}
For all $\Delta_{1},\Delta_{2} \in (s_{i}, t_{i+1}]$, 
$\Delta_{1} < \Delta_{2}$, and $y \in {\cal D}_{r}$, 
we obtain
\begin{eqnarray*}
\mathbb{E}\|({\cal E}_{1}y)(\Delta_{2}) - ({\cal E}_{1}y)(\Delta_{1})\|^{2} 
& \leq & \mathbb{E}\|({\cal T}_{q}(\Delta_{2}-s_{i})
-{\cal T}_{q}(\Delta_{1}-s_{i})){\cal K}_{i}(s_{i},g_{s_{i}}+\bar{y}_{s_{i}})\|^{2}.
\end{eqnarray*}
Since ${\cal T}_{q}$ is strongly continuous, it allows us to conclude that 
\begin{eqnarray}
\lim_{n \rightarrow \infty}\|{\cal T}_{q}(\Delta_{2}-s_{i})
-{\cal T}_{q}(\Delta_{1}-s_{i})\|^{2}=0.
\label{zz1}
\end{eqnarray}
Equations~(\ref{zz}) and (\ref{zz1}) with~(H9)(ii) imply that the operator ${\cal E}_{1}$ 
is equicontinuous on ${\cal D}_{r}$. Finally, combining steps 1--4 together 
with Ascoli's theorem, we conclude that the operator ${\cal E}_{1}$ is completely continuous.

Step 5. The operator ${\cal E}_{2}$ is a contraction map.
For $y,y^{*} \in {\cal D}_{r}$ and for $t \in (t_{i},s_{i}]$, $i=1,2,\ldots,m$, we have
\begin{eqnarray}
\mathbb{E}\|({\cal E}_{2}y)(t)-({\cal E}_{2}y^{*})(t)\|^{2}=0.\qquad\qquad
\end{eqnarray}
Similarly, for $y,y^{*} \in {\cal D}_{r}$ 
and for $t \in (s_{i},t_{i+1}]$, $i=0,1,\ldots,m$, we have
\begin{eqnarray*}
\mathbb{E}\|({\cal E}_{2}y)(t)-({\cal E}_{2}y^{*})(t)\|^{2} 
&\leq & 2\mathbb{E} \left\| \int_{s_{i}}^{t} {\cal S}_{q}
(t-e)({\cal F}(e,g_{e}+\bar{y}_{e})
-{\cal F}(e,g_{e}+\bar{y}^{*}_{e}))de\right\|^{2}\qquad\qquad\quad\\
&&+ 2\mathbb{E}\left\|\int_{s_{i}}^{t} {\cal S}_{q}(t-e)({\cal G}(e,g_{e}
+\bar{y}_{e})-{\cal G}(e,g_{e}+\bar{y}^{*}_{e}))d\hat{\cal W}(e)\right\|^{2}\\
&\leq & 2{\cal M}_{2}^{2}\varpi^{2}\left(\dfrac{N_{\cal F}T^{2q}}{q^{2}}
+\dfrac{N_{\cal G}T^{2q-1}}{2q-1}\right)\|y-y^{*}\|^{2}_{{\cal D}_{T}^{0}}.
\end{eqnarray*}
Hence,
\begin{eqnarray}
\mathbb{E}\|({\cal E}_{2}y)(t)-({\cal E}_{2}y^{*})(t)\|^{2} 
\leq 2{\cal M}_{2}^{2}\varpi^{2}\left(\dfrac{N_{\cal F}T^{2q}}{q^{2}}
+\dfrac{N_{\cal G}T^{2q-1}}{2q-1}\right)\|y-y^{*}\|^{2}_{{\cal D}_{T}^{0}}.
\qquad\qquad
\end{eqnarray}
From above, we obtain
\begin{eqnarray*}
\|{\cal E}_{2}y-{\cal E}_{2}y^{*}\|^{2}_{{\cal D}_{T}^{0}} 
\leq L_{\cal HR}\|y-y^{*}\|^{2}_{{\cal D}_{T}^{0}}.
\end{eqnarray*}
Thus, ${\cal E}_{2}$ is a contraction map. 
By Krasnoselskii's fixed point theorem, we obtain 
that problem~(\ref{eq1}) has at least one solution 
on $(-\infty, T]$.
\end{proof}

% ---------------------------------------------------------

\section{Approximate Controllability}
\label{sec:04}

We consider the following control system:
\begin{adjustwidth}{-\extralength}{0cm}
\begin{equation}
\begin{cases}
^cD_{t}^{q}z(t) = {\cal P}z(t)+{\cal A}\hat{u}(t)
+{\cal F}(t,z_{t})+{\cal G}(t,z_{t})\dfrac{d\hat{\cal W}(t)}{dt}
+\sigma(t)\dfrac{d{\cal B}^{\hat{\cal H}}(t)}{dt},
\quad t\in \cup_{i=0}^{m}(s_{i},t_{i+1}],& \text{}\\
z(t)  =  {\cal K}_{i}(t,z_{t}),
\quad t\in \cup_{i=1}^{m}(t_{i},s_{i}], &\text{}\\
z(t)= \phi(t),\quad \phi(t)\in {\cal D}_{h}. & \text{}\label{eq2}
\end{cases}
\end{equation}
\end{adjustwidth}
The control $\hat{u}(\cdot)\in L^{2}({\cal J},{\cal U})$, 
where $L^{2}({\cal J},{\cal U})$ is the Hilbert space 
of all admissible control functions. The operator ${\cal A}$ is linear and bounded 
from the separable Hilbert space ${\cal U}$ into ${\cal Z}$. 
Assume that the linear system 
\begin{equation}
\begin{cases}
^cD_{t}^{q}z(t) = {\cal P}z(t)+{\cal A}\hat{u}(t),
\quad t\in [0,T],& \text{}\\
z(t)= \phi(t),
\quad \phi(t)\in {\cal D}_{h}. & \text{}\label{eq3}
\end{cases}
\end{equation}
Define the operator $t_{s_{i}}^{t_{i+1}}$ 
associated with system of~(\ref{eq3}) as
\begin{eqnarray*}
{t _{s_{i}}^{t_{i+1}} = \int_{s_{i}}^{t_{i+1}} 
{\cal S}_{q}(t_{i+1}-e){\cal A}{\cal A}^{*}{\cal S}_{q}^{*}(t_{i+1}-e)de}.
\end{eqnarray*} 
Here, ${\cal A}^{*}$ and ${{\cal S}_{q}^{*}(t)}$ are the adjoint of 
${\cal A}$ and ${\cal S}_{q}(t)$, respectively. The
operator $t _{s_{i}}^{t_{i+1}}$ is a bounded and linear operator. 

\begin{Definition} 
System~(\ref{eq2}) is approximately controllable on $[0,T]$ if  
$\overline{{\cal R}(T,\phi,\hat{u})}=L^{2}({{\cal F}_{T},{\cal Z}})$, 
where ${\cal R}(T,\phi,\hat{u})=\{z(\phi,\hat{u})(T):z\; 
\mbox{is the solution of problem}~(\ref{eq2})\; \mbox{and}\; 
\hat{u}\in L^2({\cal J},{\cal U}) \}$.
\end{Definition}
 
The following assumption is needed.
\begin{enumerate}
\item [{[AC]:}]~System~(\ref{eq3}) is approximate controllability on ${\cal J}$.
\end{enumerate}

Note that system~(\ref{eq3}) is approximately controllable on ${\cal J}$ only if 
\begin{eqnarray}
{\Delta(\Lambda, t _{s_{i}}^{t_{i+1}}) 
= (\Lambda I+t _{s_{i}}^{t_{i+1}})^{-1}}
\rightarrow 0\;\mbox{as}\;\Lambda \rightarrow 0.
\label{AC}
\end{eqnarray}
\begin{Definition} 
An ${\cal F}_{t}$-adapted random process $z:(-\infty,T] \rightarrow {\cal Z}$ 
is called the mild solution of~(\ref{eq2}) if for every $t \in {\cal J}$, 
$z(t)$ satisfies $z_{0}=\phi\in {\cal D}_{h},$ $z(t)= {\cal K}_{i}(t,z_{t})$ 
for all $t \in (t_{i},s_{i}]$, $i=1,2,\ldots,m$, and 
\begin{eqnarray*}
z(t)&=&  \int_{0}^{t} {\cal S}_{q}(t-e)[{\cal F}(e,z_{e})+{\cal A}\hat{u}(e)]de\\
&&+\int_{0}^{t}{\cal S}_{q}(t-e){\cal G}(e,z_{e})d\hat{\cal W}(e)
+ \int_{0}^{t}{\cal S}_{q}(t-e)\sigma(e)d{\cal B}^{\hat{\cal H}}(e),
\end{eqnarray*}
for all $t \in [0,t_{1}]$, and
\begin{eqnarray}
z(t)&=& {\cal T}_{q}(t-s_{i}){\cal K}_{i}(s_{i},z_{s_{i}}) 
+ \int_{s_{i}}^{t} {\cal S}_{q}(t-e)[{\cal F}(e,z_{e})+{\cal A}\hat{u}(e)]de\nonumber\\
&&+\int_{s_{i}}^{t} {\cal S}_{q}(t-e){\cal G}(e,z_{e})d\hat{\cal W}(e) 
+ \int_{s_{i}}^{t} {\cal S}_{q}(t-e)\sigma(e)d{\cal B}^{\hat{\cal H}}(e),\label{mild2}
\end{eqnarray}
for all $t \in (s_{i},t_{i+1}]$, $i=1,2,\ldots,m$.
\end{Definition}

\begin{Lemma} 
For any $z_{t_{i+1}}$ $\in$  $L^{2}({\cal F}_{T},{\cal Z})$, 
there exist $\phi_{1} \in L^{2}(\Omega,L^{2}([s_{i},t_{i+1}],
{\cal L}_{2}^{1}({\cal Y}_{1},{\cal Z}))$ and 
$\phi_{2} \in L^{2}([s_{i},t_{i+1}],{\cal L}^{2}_{2}({\cal Y}_{2},{\cal Z}))$ 
such that
\begin{eqnarray*}
z_{t_{i+1}}= \mathbb{E}z_{t_{i+1}}
+\int_{s_{i}}^{t_{i+1}}\phi_{1}(e)d\hat{\cal W}(e)
+\int_{s_{i}}^{t_{i+1}}\phi_{2}(e)d{\cal B}^{\hat{\cal H}}(e).
\end{eqnarray*}
\end{Lemma}

Next, we choose the control $\hat{u}^{\Lambda}(t)$ as follows:
\begin{eqnarray}
\hat{u}^{\Lambda}(t) = {\cal A}^{*}{\cal S}_{q}^{*}(t_{i+1}-t)
\Delta(\Lambda, t _{s_{i}}^{t_{i+1}}) p(z(\cdot)),
\end{eqnarray}
where
\begin{eqnarray*}
p(z(\cdot))&=& z_{t_{i+1}}-{\cal T}_{q}(t_{i+1}-s_{i}){\cal K}_{i}(s_{i},z_{s_{i}}) 
- \int_{s_{i}}^{t_{i+1}} {\cal S}_{q}(t_{i+1}-e){\cal F}(e,z_{e})de\nonumber\\
&&-\int_{s_{i}}^{t_{i+1}} {\cal S}_{q}(t_{i+1}-e){\cal G}(e,z_{e})d\hat{\cal W}(e) 
- \int_{s_{i}}^{t_{i+1}} {\cal S}_{q}(t_{i+1}-e)\sigma(e)d{\cal B}^{\hat{\cal H}}(e),\\
&&\forall \; t \in (s_{i},t_{i+1}],\quad i=0,1,\ldots,m,
\end{eqnarray*}
and ${\cal K}_{0}(0,\cdot)=0$, $z(t_{m+1}) = z_{t_{m+1}}=z_{T}$.

\begin{Theorem}
Assume the hypotheses {(H1)--(H9)} are satisfied. 
Then, the problem~(\ref{eq2}) has at least one mild solution 
on $(-\infty, T]$.
\end{Theorem}

\begin{proof}
The proof is a consequence of Theorem~\ref{xxx}.
\end{proof}

\begin{Theorem} 
Assume that the hypotheses {(H1)--(H9)} and [AC] are satisfied. 
Then functions ${\cal F}$ and ${\cal G}$ are uniformly bounded on 
their respective domains. Moreover, the system~(\ref{eq2}) 
is approximately controllable on $[0,T]$.
\end{Theorem}

\begin{proof}
Let $z^{\Lambda}$ be a fixed point of ${\cal E}_{1}+{\cal E}_{2}$. 
Using Fubini's theorem, we get
\begin{eqnarray}
z^{\Lambda}(t_{i+1}) = z_{t_{i+1}}- \Lambda 
\Delta(\Lambda, t _{s_{i}}^{t_{i+1}})p(z^{\Lambda}(\cdot)),
\label{22}
\end{eqnarray}
where
\begin{eqnarray*}
p(z^{\Lambda}(\cdot))
&=& z_{t_{i+1}}-{\cal T}_{q}(t_{i+1}-s_{i}){\cal K}_{i}(s_{i},z^{\Lambda}_{s_{i}}) 
- \int_{s_{i}}^{t_{i+1}} {\cal S}_{q}(t_{i+1}-e){\cal F}(e,z^{\Lambda}_{e})de\nonumber\\
&&-\int_{s_{i}}^{t_{i+1}} {\cal S}_{q}(t_{i+1}-e){\cal G}(e,z^{\Lambda}_{e})d\hat{\cal W}(e)
-\int_{s_{i}}^{t_{i+1}} {\cal S}_{q}(t_{i+1}-e)\sigma(e)d{\cal B}^{\hat{\cal H}}(e),\\
&&\forall \; t \in (s_{i},t_{i+1}],\;i=0,1,\ldots,m.
\end{eqnarray*}
The functions ${\cal F}$ and ${\cal G}$ are uniformly bounded. Hence, there exists 
a subsequence, still represented by ${{\cal F}(e,z^{\Lambda}_{e})}$ 
and ${{\cal G}(e, z^{\Lambda}_{e})}$, that weakly converge to, say, 
${\cal F}(e)$ and ${\cal G}(e)$ in ${\cal Z}$ 
and ${\cal L}^{1}_{2}({\cal Y}_{1},{\cal Z})$, respectively. Let us define
\begin{eqnarray*}
\eta &=& z_{t_{i+1}} - {\cal T}_{q}(t_{i+1}-s_{i}){\cal K}_{i}(s_{i},z_{s_{i}}) 
- \int_{s_{i}}^{t_{i+1}} {\cal S}_{q}(t_{i+1}-e){\cal F}(e)de\nonumber\\
&&-\int_{s_{i}}^{t_{i+1}}{\cal S}_{q}(t_{i+1}-e){\cal G}(e)d\hat{\cal W}(e)
- \int_{s_{i}}^{t_{i+1}} {\cal S}_{q}(t_{i+1}-e)\sigma(e)d{\cal B}^{\hat{\cal H}}(e),\\
&&\forall \; t \in (s_{i},t_{i+1}],\;i=0,1,\ldots,m.
\end{eqnarray*}
For $t \in (s_{i},t_{i+1}]$, $i=0,1,\ldots,m$, we have
\begin{eqnarray*}
\mathbb{E}\| p(z^{\Lambda})-\eta \|^{2} 
&\leq & 3 \mathbb{E}\|{\cal T}_{q}(t_{i+1}-s_{i})({\cal K}_{i}(s_{i},z^{\Lambda}_{s_{i}})
-{\cal K}_{i}(s_{i},z_{s_{i}}))\|^{2}\nonumber \\&&
+ 3\mathbb{E}\left\|\int_{s_{i}}^{t_{i+1}} {\cal S}_{q}(t_{i+1}-e)({\cal F}(e,z^{\Lambda}_{e})
-{\cal F}(e))de\right\|^{2}\nonumber\\
&&+3\mathbb{E}\left\|\int_{s_{i}}^{t_{i+1}} {\cal S}_{q}(t_{i+1}-e)
({\cal G}(e,z^{\Lambda}_{e})-{\cal G}(e))d\hat{\cal W}(e)\right\|^{2}.
\end{eqnarray*}
By the infinite dimensional version of the Arzela--Ascoli theorem, 
we obtain that 
$$
\bar{k}(\cdot) \rightarrow \int_{\cdot}^{\cdot} {\cal S}_{q}(\cdot-e)\bar{k}(e)de
$$ 
is a compact operator. For all $t \in [0,T]$, 
\begin{eqnarray}
\label{es}
\mathbb{E}\| p(z^{\Lambda})-\eta \|^{2} 
\rightarrow 0\; \mbox{as}\;\Lambda \rightarrow 0^{+}.
\end{eqnarray}
By Equation~(\ref{22}), we get
\begin{eqnarray*}
\mathbb{E}\| z^{\Lambda}(t_{i+1}) - z_{t_{i+1}} \|^{2} 
\leq  \mathbb{E}\|\Lambda \Delta(\Lambda, t _{s_{i}}^{t_{i+1}})(\eta)\|^{2}
+\mathbb{E}\|\Lambda \Delta(\Lambda, t _{s_{i}}^{t_{i+1}})\|^{2}
\mathbb{E}\| p(z^{\Lambda})-\eta \|^{2}.
\end{eqnarray*}
By (\ref{AC}) and (\ref{es}), we get 
$$
\mathbb{E}\|z^{\Lambda}(t_{i+1}) - z_{t_{i+1}} \|^{2} 
\rightarrow 0\; \mbox{as}\;\theta \rightarrow 0^{+}.
$$ 
Thus, the system~(\ref{eq2}) is approximate 
controllable on the interval $[0,T]$.
\end{proof}

% ---------------------------------------------------------

\section{Example}
\label{sec:05}

We consider the following fractional stochastic control system:
\begin{equation}
\begin{cases}
^cD_{t}^{q} y(t,z)=\dfrac{\partial^{2}}{\partial z^{2}} y(t,z) + \Theta(t, z)
+\int_{-\infty}^{t}e^{4(r-t)}y(r,z)dr\\\qquad\qquad\qquad
+\int_{-\infty}^{t}e^{6(r-t)}y(r,z)dr\dfrac{d\hat{\cal W}(t)}{dt}
+P(t) \dfrac{d{\cal B}^{\hat{\cal H}}(t)}{dt},\\
\qquad\qquad\qquad\;\; y \in (0,\pi),\;t\in[2i,2i+1],
\;i=0,1,\ldots,m, &\text{}\\
y(t,z) = \int_{-\infty}^{t}G_{i}(r-t)y(r,z)dr,\;t\in(2i-1,2i],
\;i=1,2,\ldots,m, &\text{}\\
y(t,0)=\;0\;=y(t,\pi), & \text{}\\
y(t,z)= \phi(t,z),\;t \in (-\infty,0], & \text{}\label{E}
\end{cases}
\end{equation}
where $^cD_{t}^{q}$ is the Caputo derivative of order $1/2 < q < 1$, 
$0 = s_0 = t_0 < t_1 < s_1 < t_2<\cdots <t_m < s_m < t_{m+1} 
= T <\infty$ with $s_{i}=2i,\;t_{i}=2i-1$.  

Let ${\cal Z}=L^{2}([0,\pi])$ and the operator ${\cal P}$ be defined by
\begin{eqnarray*}
{\cal P}w=w''\;,\;\;{\cal D}({\cal P})=H^{2}(0,\pi)\cap H^{1}_{0}(0,\pi).
\end{eqnarray*}
Clearly, ${\cal P}$ is the generator of an analytic semigroup 
$\{{\cal S}(t):t \geqslant 0\}$. The spectral representation of ${\cal S}(t)$ 
is given by
$$ 
{\cal S}(t)w=\sum_{n \in \mathbb{N}}e^{-n^{2}t}
\langle w,w_{n}\rangle w_{n},
$$
where
\begin{eqnarray*}
w_{n}(y)= \sqrt{2/\pi}\, \sin\,(ny),
\; n \in \mathbb{N},
\end{eqnarray*}
is the orthogonal set of eigenvectors corresponding to the eigenvalue 
$\lambda_{n}=-n^{2}$ of ${\cal P}$. The semigroup $\{{\cal S}(t):t\geq 0\}$ 
is compact and uniformly bounded, so that ${\cal R}(\lambda,{\cal P}) 
= (\lambda I-{\cal P})^{-1}$ is a compact operator for all $\lambda \in \rho({\cal P})$, 
i.e., ${\cal P} \in {\cal P}^{q}(\theta_{0}, \omega_{0})$. Let $h(e) = e^{2e}$, $e < 0$. 
Then $\varpi = \int_{-\infty}^{0} h(e)de = 1/2$ and we define
\begin{eqnarray*}
\|\phi\|_{{\cal D}_{h}}= \int_{-\infty}^{0} h(e)
\sup_{e\leq\theta\leq 0}(\mathbb{E}|\phi(\theta)|^{2})^{1/2}de,
\quad \phi\in {\cal D}_{h}.
\end{eqnarray*}
Hence, $(t,\phi) \in [0, T] \times {\cal D}_{h}$. The bounded linear 
operator ${\cal A}$ is defined by ${\cal A}\hat{u}(t)(z)=\Theta(t, z)$.

Define the functions ${\cal F} : {\cal J} \times {\cal D}_{h} 
\rightarrow {\cal Z}$, ${\cal G} : {\cal J} \times {\cal D}_{h} 
\rightarrow L_{2}({\cal Y}_{1},{\cal Z})$, and 
${\cal K}_{i} : (t_{i},s_{i}] \times {\cal D}_{h} \rightarrow {\cal Z}$ as 
\begin{eqnarray*}
{\cal F}(t,\phi)(z) &=&\int_{-\infty}^{0}e^{4\theta}(\phi(\theta)(z))d\theta,\\
{\cal G}(t,\phi)(z)&=&\int_{-\infty}^{0}e^{6\theta}(\phi(\theta)(z))d\theta,\\
{\cal K}_{i}(t,\phi)(z)&=& \int_{-\infty}^{0}G_{}(\theta)(\phi(\theta)(z))d\theta.
\end{eqnarray*}
Assume that 
$$
\int_{0}^{T}\|\sigma(e)\|^{2}_{{\cal L}^{2}_{2}}de < \infty.
$$  
The system~(\ref{E}) can be written as an abstract formulation of (\ref{eq1}), 
and thus previous theorems can be applied to guarantee both existence 
and approximate controllability results.

% ---------------------------------------------------------

\section{Conclusions}
\label{sec:06}

We have investigated impulsive fractional stochastic control systems
defined on separable Hilbert spaces. The proposed problem is driven 
by mixed noise, i.e., it involves both a $Q$-Wiener process 
and a $Q$-fractional Brownian motion with the Hurst parameter 
$\hat{\cal H} \in (1/2, 1)$. For our results, we have mainly applied 
fixed point techniques, a $q$-resolvent family, and fractional calculus. 
The obtained results are supported by an illustrative example. 
As further directions of investigation and continuation to this work,
it would be interesting to investigate the sensitivity on the noise range
and develop numerical and computational methods to approximate
the solution. We also intend to extend our results 
via discrete fractional calculus. 

% ---------------------------------------------------------

\vspace{6pt} 

% ---------------------------------------------------------

\authorcontributions{Conceptualization, N.H. and R.D.; 
methodology, R.D.; validation, A.D.; 
formal analysis, D.F.M.T.; investigation, A.D.; 
writing---original draft preparation, N.H. and R.D.; 
writing---review and editing, A.D. and D.F.M.T.; 
supervision, A.D.; project administration, D.F.M.T. 
All authors have read and agreed to the published version of the manuscript.}

% ---------------------------------------------------------

\funding{Debbouche and Torres were supported by 
\emph{Funda\c{c}\~{a}o para a Ci\^{e}ncia e a Tecnologia} (FCT) 
within project number UIDB/04106/2020 (CIDMA).} 

% ---------------------------------------------------------

\acknowledgments{The authors are grateful to three reviewers 
for their constructive comments and suggestions.} 

% ---------------------------------------------------------

\conflictsofinterest{The authors declare no conflict of interest. 
The funder had no role in the design of the study; in the collection, analyses, 
or interpretation of data; in the writing of the manuscript; 
or in the decision to publish the results.}

% ---------------------------------------------------------

\begin{adjustwidth}{-\extralength}{0cm}

\reftitle{References}

% ---------------------------------------------------------

\PublishersNote{}
\end{adjustwidth}
\end{document}